\documentclass[10pt,draft,a4paper,reqno]{amsart}
\usepackage[latin1]{inputenc}
\usepackage{amsmath,amssymb,amscd,amsthm,amsxtra,verbatim,ifthen}
\usepackage[english]{babel}
\usepackage[all]{xy}

\mathchardef\ordinarycolon\mathcode`\: \mathcode`\:=\string"8000
\begingroup \catcode`\:=\active
  \gdef:{\mathrel{\mathop\ordinarycolon}}
\endgroup

\theoremstyle{plain}

\newtheorem*{singletheorem}{Theorem} % unnumeriert

\theoremstyle{definition}

\theoremstyle{remark}

\newtheorem*{singleremark}{Remark}
\newtheorem{example}{Example}

\newcommand{\im}{\mathrm{im}}
\newcommand{\pr}{\mathrm{pr}}
\newcommand{\id}{\mathrm{id}}

\newcommand{\reg}{\mathrm{reg}}

\newcommand{\Ad}{\mathrm{Ad}}
\newcommand{\ad}{\mathrm{ad}}
\newcommand{\vol}{\mathrm{vol}}
\newcommand{\iso}{\mathrm{Iso}}

\newcommand{\R}{\mathbf{R}}
\newcommand{\C}{\mathbf{C}}

\newcommand{\N}{\mathbf{N}}

\newcommand{\kc}{{\mathcal C}}
\newcommand{\ku}{{\mathcal U}}

\newcommand{\gotha}{{\mathfrak a}}
\newcommand{\gothg}{{\mathfrak g}}
\newcommand{\gothh}{{\mathfrak h}}
\newcommand{\gothk}{{\mathfrak k}}
\newcommand{\gothl}{{\mathfrak l}}
\newcommand{\gothm}{{\mathfrak m}}

\newcommand{\gothp}{{\mathfrak p}}

\newcommand{\gotht}{{\mathfrak t}}

\newcommand{\gothz}{{\mathfrak z}}

\begin{document}
\title{An Integration Formula for Polar Actions}
\author{Frederick Magata}
\begin{abstract}
We prove an analogue of Weyl's Integration Formula for compact Lie
groups in the context of polar actions. We also show how certain
classical examples from the literature can be viewed as special
cases of our result.
\end{abstract}
\maketitle

\section{The integration formula}
Let $\varphi:G\times M\to M$ be an smooth and proper isometric
action action of a Lie group $G$ on a Riemannian
manifold $M$. % In particular, the action $\varphi$ is effective.
$\varphi$ is called \emph{polar}, if there exists an embedded
submanifold $\Sigma$ of $M$ which intersects with every $G$-orbit
and is orthogonal to each orbit in the intersection points. Our main
result is:

\begin{singletheorem} Let $\varphi$ be a polar action and
$\Sigma\subset M$ a section. We put $H:=Z_G(\Sigma)$ and denote by
$W=W(\Sigma)$ the generalized Weyl group of $\Sigma$. Furthermore,
let $\omega_s:G/H\to M,\, gH\mapsto g\cdot s$ denote the \lq orbit
map\rq\ in the point $s\in M$ and assume $W$ to be finite.
%By abuse of notation, we denote with $\varphi$ also the induced
%map from $G/H\times\Sigma$ to $M$.
Then the following holds:

\begin{enumerate}

\item The function $\delta(s):=|\det(d\omega_s(eH))|$ on $\Sigma$
is continuous and $W$-in\-va\-ri\-ant. If $G$ is compact, then
$\delta$ is a volume scaling function in the following way:
$$\delta(s)=\left\{\begin{array}{ll}
0 & \text{if } s \text{ is singular,}\\
|G_s/H|\cdot\frac{\vol_l(G\cdot s)}{\vol_l(G/H)}\cdot  & \text{if }
s \text{ is regular or exceptional,}\end{array}\right.
$$
where $l:=\dim(G/H)$ is the dimension of a principal orbit.

\item The assignment $\Psi:f\in\kc_c(M)\mapsto F\in\kc_c(G/H\times
\Sigma)^W$, given by $$F(gH,s):=f(g\cdot s)\, \delta(s),$$ extends
to a continuous isomorphism from $L^1(M)$ onto
$L^1(G/H\times\Sigma)^W$.

\item For any $f\in L^1(M)$ we have the formula
$$\int_M f(x)\, dx=\frac{1}{|W|}\int_\Sigma \left(\int_{G/H} f(g\cdot
s)\, d(gH)\right)\delta(s)\, ds,$$
where $d(gH)=dg/dh$ is induced by the Haar measures $dg$ on $G$, resp. $dh$ on
$H$. It coincides with the Riemannian measure of $G/H$ with respect to the
quotient metric.

\item Let $G$ be compact and $c:=\vol_l(G/H)$. Then $f\mapsto
\sqrt[p]{\frac{c\delta}{|W|}}f|_\Sigma$ is an isometry from
$L^p(M)^G$ to $L^p(\Sigma)^W$ and for any $f\in L^1(M)^G$ we have
the formula $$\int_M f(x)\, dx=\frac{c}{|W|}\int_\Sigma
f(s)\delta(s)\,ds.$$

\end{enumerate}
\end{singletheorem}

In particular, (iii) is a generalization of Weyl's celebrated
integration formula.

Now some remarks concerning our notation are in order:
\begin{singleremark}
{\ }\newline \vspace{-3ex}
\begin{enumerate}
\item Unless otherwise stated, we always consider the Riemannian measure
associated to the Riemannian metric on a Riemannian manifold. In
particular, we do not necessarily assume our manifolds to be
oriented.

\item We equip Lie groups with a left invariant Riemannian metric
and the corresponding Haar measure. The latter one then coincides
with the Riemannian measure. For a compact Lie group, we further
assume its invariant metric normalized such that its total volume
becomes one.

\item For function spaces on manifolds with some $G$-action, a
superscript denotes the subset of $G$-invariant functions. E.g.
$\kc(M)^G$ denotes the set of all $G$-invariant continuous
functions on $M$.

\item The presentation of the theorem follows Theorem (3.14.1) and
Corollary (3.14.2) of \cite{Duistermaat&Kolk}, whose proofs include
most of the necessary ideas for the generalized version we present
in this paper. In \cite{Helgason2}(Ch. I, \S5), several important
special cases of the integral formul\ae\ are proved and explicit
formulas for the Jacobians $\delta$ are given there. We have
included some of them in section \ref{s:examples} below.

\item For the proof of the theorem, we expect the reader to be
familiar with certain basic properties of polar actions and
sections which can be found in \cite{Palais&Terng}, or in the more
recent literature \cite{Berndt}. Here are some of them:
\begin{enumerate}
\item The set of $G$-regular points lies open and dense in every
section.

\item Sections are totally geodesic

\item The slice representation in each point is polar.

\item The generalized Weyl group
$W(\Sigma)=N_G(\Sigma)/Z_G(\Sigma)$ parameterizes the intersection
of every $G$-orbit with the section $\Sigma$.

\item The generalized Weyl group is a discrete group. The
assumption on the finiteness of $W$ in the theorem is met in the
cases that $G$ is a compact group, or if $\Sigma$ is compact.

\item The $G$-regular set $M^\reg$ of $M$ is covered $|W|$-fold by
$G/H\times\Sigma^\reg$ where $H=Z_G(\Sigma)$. In case of a compact
Lie group acting on itself via conjugation, this fact is sometimes
called \lq Weyl's covering theorem\rq.

\end{enumerate}

\item The assumption that $G$ acts properly, smoothly and in an
isometric fashion on $M$ means that we have a smooth homomorphism
$G\to\iso(M)$ whose image is a closed subgroup of the isometry group
of $M$. In particular, the action need not be effective as it is the
case with an arbitrary compact Lie group acting on itself via
conjugation (see the section \lq Examples\rq\ below).

\end{enumerate}
\end{singleremark}

\begin{proof}
To (i): We first deal with the invariance properties of $\delta$.
For this purpose, and for the integration formula in (iii), we
need the following calculations. By abuse of notation, we denote
by $\varphi$ also the induced map $G/H\times \Sigma\to M$. For
general $g\in G$, the defining properties of a group action imply
\begin{equation}
\varphi\circ(l_g\times\id_\Sigma)=\phi_g\circ\varphi,
\end{equation}
where $l_g:G/H\to G/H$ is the left translation by $g$ and
similarly $\phi_g:M\to M$ is the isometry of $M$ induced by $g$
via the group action $\varphi$. Forming the differential in
$(eH,s)$, we obtain
\begin{equation}
d\varphi(gH,s)\circ(dl_g(eH)\times\id_{T_s\Sigma})=d\phi_g(s)\circ
d\varphi(eH,s)
\end{equation}
which we rewrite in the form
$$d\varphi(gH,s)=d\phi_g(s)\circ d\varphi(eH,s)\circ(dl_{g^{-1}}(gH)\times
\id_{T_s\Sigma}).$$

If in (1) we have that $g\cdot s$ lies in $\Sigma$ again, we may
further assume that $g\in N_G(\Sigma)$, because the intersection
$\Sigma\cap (G\cdot s)$ is parameterized by the normalizer
$N_G(\Sigma)$ of $\Sigma$. For such $g\in N_G(\Sigma)$, we get
another equation similar to (2) in the following way. We first
observe that:
$$\varphi\circ(r_g\times\id_\Sigma)=\varphi\circ(\id_{G/H}\times\phi_g),$$
where $r_g:G/H\to G/H,\ kH\mapsto kgH$, the right translation, is
well defined since $g\in N_G(\Sigma)$. Forming the derivative in
$(eH,s)$ again, yields
\begin{equation}
d\varphi(gH,s)\circ(d
r_g(eH)\times\id_{T_s\Sigma})=d\varphi(eH,g\cdot
s)\circ(\id_{\gothg/\gothh}\times d\phi_g(s)). \end{equation}
Combining (2) and (3), we obtain the following formula for
$d\varphi(eH,g\cdot s)$:
\begin{equation}
d\varphi(eH,g\cdot s)=d\phi_g(s)\circ
d\varphi(eH,s)\circ(\Ad_{g^{-1}}\times d\phi_g^{-1}(s))
\end{equation}
In other words, the differentials of $\varphi$ in $(eH,s)$ and
$(eH,g\cdot s)$ for $g\in N_G(\Sigma)$ differ only by the
isometries $d\phi_g(s)$ of $M$ and $(\Ad_{g^{-1}}\times
d\phi_g^{-1}(s))$ of $G/H\times\Sigma$. Note that the latter one
is true, because $g\in N_G(\Sigma)$ and left translation by $g$ is
an isometry of $G$ as well as $G/H$ as a consequence of the remark
(ii) above.

Recall that $\det(d\varphi(gH,s))$ is computed by choosing some
orthonormal bases on $T_{gH}G/H\times T_s\Sigma$ and $T_{g\cdot s}
M$ and then computing the usual determinant of the matrix
representation of $d\varphi(gH,s)$ with respect to the given
bases. The determinant obtained in this way is, up to a sign,
independent of the choice of the bases. Now,
$$d\varphi(eH,s)(X+\gothh,v)=d\omega_s(eH)(X)+d\phi_e(s)(v)=d\omega_s(eH)(X)+v.$$
By choosing orthonormal bases $x_1,\dots,x_l$ of $\gothg/\gothh$
and $y_{l+1},\dots, y_{l+m}$ of $T_s\Sigma$ and completing the
latter to an on-basis {$y_1,\dots,y_l,y_{l+1},\dots y_{l+m}$} on
$T_s M$, we obtain the following matrix representation of
$d\varphi(eH,s)$ with respect to these bases:
$$\left(\begin{array}{cc}
d\omega_s(eH)&0\\
0 & \id_{T_s\Sigma}
\end{array}\right).$$

We therefore have
$\delta(s)=|\det(d\omega_s(eH))|=|\det(d\varphi(eH,s))|$, which we
will need in the proof of the integral formula (iii); and by (4),
we have that $\delta(s)=\delta(g\cdot s)$ for any $g\in
N_G(\Sigma)$; proving the $W$-invariance of $\delta$.

In order to show the continuity of $\delta$, we express $\delta$
as a composition of continuous mappings by means of local frames.
More precisely, consider an on-basis $x_1,\dots, x_l$ on
$\gothg/\gothh$ and a local frame $Y_1,\dots,
Y_l,Y_{l+1},\dots,Y_{m+l}$ of $M$ on $U\subset\Sigma$ about $s$
which is adapted to $\Sigma$; i.e. $Y_{l+1},\dots, Y_{l+m}$ is
tangent to $\Sigma$. Mapping $x_i$ to $Y_i(p)$, we obtain for
every $p\in U$ a linear isometry $L_p:\gothg/\gothh\to
\nu_p\Sigma$. This yields a smooth mapping
$L:U\times\gothg/\gothh\to \nu(U)$ which is nothing but a local
trivialization of the normal bundle $\nu(\Sigma)$ of $\Sigma$.
Moreover, the pullback via $L_p$ of the measure on $\nu_p\Sigma$
induced by $dx_p$ gives $d(gH)_{eH}$. By means of $L$, we may
write $\delta(p)=|\det(L_p^{-1}\circ d\omega_p(eH))|$ for all
$p\in U$ where $\det(\cdot)$ now denotes the usual determinant for
endomorphisms of $\gothg/\gothh$. This proves the continuity of
$\delta$. In fact, $L_p$ and $d\omega_p(eH)$ depend smoothly on
$p$ on $U$; the latter because $\varphi$ is a smooth action and
$d\omega_p(eH)=\frac{\partial}{\partial
p}\varphi(eH,p)|_{\gothg/\gothh\times\{0\}}$.

Now if $s$ is $G$-regular, then $\omega_s$ is an embedding, if $s$
is exceptional, $\omega_s$ is a $|G_p/H|$-fold covering and if $s$
is singular, $\omega_s$ is a bundle over $G\cdot s$ with fibre
$G_s/H$. If $s$ is singular, then $l>\dim\, G\cdot s=\dim\,
\im(d\omega_s(eH))$ and hence the rank of $d\omega_s(eH)$ is
strictly less than $l$. Thus $\delta(s)=0$ in that case.

If $s$ is $G$-regular, then by definition $\vol_l(G\cdot
s)=\int_{G\cdot s} 1\, dy$ where $dy$ denotes the Riemannian measure
on $G\cdot s$ related to the metric on $G\cdot s$ induced from $M$.
Applying the transformation theorem to the embedding
$\omega_s:G/H\to G\cdot s$, we have:
%\begin{eqnarray*}
$$\vol_l(G\cdot s)=\int_{G\cdot s}\hspace{-1em}1\,
dy=\int_{G/H}\hspace{-1em}|\det(d\omega_s(gH))|\, d(gH)
=\int_{G/H}\hspace{-1em}\delta(s)\, d(gH)=\vol_l(G/H)\delta(s).$$
%\end{eqnarray*}
A similar consideration for exceptional $s\in\Sigma$ yields:
$$\vol_l(G\cdot s)=\frac{1}{|G_s/H|}\int_{G/H}\hspace{-1em}|\det(d\omega_s(gH))|\, d(gH)=\frac{\vol_l(G/H)\delta(s)}{|G_s/H|}.$$
this concludes the proof of (i). Note that we did not use the
finiteness of $W$ here.

To (ii): First note that $\Psi$ is well defined. In fact, $F$ is
certainly well defined, continuous with compact support and for
every $(gH,s)\in G/H\times \Sigma$ and $n\in N_G(\Sigma)$ we have
$$F(gn^{-1}H,n\cdot s)=f(g\cdot s)\delta(n\cdot s)=f(g\cdot s)\delta(s).$$
That is, we used the Weyl group invariance of $\delta$ from (i).
Thus $F$ is Weyl-invariant too. It is also obvious that $\Psi$ is
a linear mapping. Concerning the continuity of $\Psi$, we should
first note that we consider $\kc_c(M)$ and
$\kc_c(G/H\times\Sigma)$ equipped with their particular
$L^1$-norm. Now let $f_n$ be a sequence of functions in $\kc_c(M)$
converging to zero. Then for $F_n=\Psi(f_n)$ we have:
\begin{eqnarray*}
\|F_n\|_1&=&\int_{G/H\times\Sigma}|F_n(gH,s)|\,
(dgH,s)=\int_\Sigma\int_{G/H}|f_n(g\cdot s)|\delta(s)\, dgH\, ds\\
&=&|W|\cdot\int_M|f_n(x)|\, dx=|W|\cdot\|f_n\|_1\to 0, \text{ as }
n\to\infty.
\end{eqnarray*}
That is, we read formula (iii) backwards. Since $\kc_c(\cdot)$ is
a dense subspace of $L^1(\cdot)$, \mbox{$\Psi$ induces a}
continuous map $L^1(M)\to L^1(G/H\times\Sigma)^W$ which again we
denote by $\Psi$. Concerning the surjectivity of $\Psi$, we
suppose that $f_1,f_2\in L^1(M)$ with $\Psi(f_1)=\Psi(f_2)$. Since
$G/H\times \Sigma^\reg$ is open and dense in $G/H\times \Sigma$,
we have for almost every $(gH,s)\in G/H\times \Sigma$:
\begin{equation}\label{psi_inj}
f_1(g\cdot s)\delta(s)=f_2(g\cdot s)\delta(s).
\end{equation}
By (i), we have even in the case that $G$ is not compact that
$\delta(s)\ne 0$ for $s\in\Sigma^\reg$. Therefore, we may divide
by $\delta(s)$ in (\ref{psi_inj}) and since
$\varphi(G/H\times\Sigma^\reg)=M^\reg$ we obtain that $f_1=f_2$
except for a set with measure zero.

Concerning the surjectivity of $\Psi$, we first define for every
$F\in\kc_c(G/H\times\Sigma^\reg)^W$ a function $f\in\kc_c(M^\reg)$
by
$$f:=\frac{F\circ\varphi^{-1}}{\delta\circ\pr_2\circ\varphi^{-1}}$$
where $\pr_2:G/H\times\Sigma^\reg\to\Sigma^\reg$ is the projection
onto the second factor. Note that the functions
$F\circ\varphi^{-1}$ and $\delta\circ\pr_2\circ\varphi^{-1}$ are
well defined since $F$ and $\delta\circ\pr_2$ are $W$-invariant
and the fibres of $\varphi$ are precisely the $W$-orbits on
$G/H\times\Sigma$. Continuity can easily be verified by
considering local trivialization for the covering $\varphi$.
Finally, $F\circ\varphi^{-1}$ has compact support as the fibres of
$\varphi$ are compact. Furthermore, $\delta$ has no zeros on
$\Sigma^\reg$ so that $\delta\circ\pr_2\circ\varphi^{-1}$ has no
zeros on $M^\reg$. In conclusion, $f$ is well defined.

We obviously have
$\Psi(f)=F$ and by (iii) we obtain: \begin{eqnarray*}
\int_{M^\reg} f(x)\,dx&=&\frac{1}{|W|}\int_{\Sigma^\reg} \left(\int_{G/H}
\frac{F(gH,s)}{\delta(s)}\, d(gH)\right)\delta(s)\,ds\\
&=&\frac{1}{|W|}\int_{\Sigma^\reg} \left(\int_{G/H} F(gH,s)\, d(gH)\right)\,ds.
\end{eqnarray*}

Now consider an arbitrary $F\in L^1(G/H\times\Sigma^\reg)^W$. We
can approximate $F$ by a sequence
$(F_n)_{n\in\N}\subset\kc_c(G/H\times\Sigma^\reg)^W$ and, as
indicated above, form a sequence
$(f_n)_{n\in\N}\subset\kc_c(M^\reg)$ with $\Psi(f_n)=F_n$. Now the
$F_n$ form a Cauchy sequence and applying the integral identity
above to $|f_n-f_m|$, we see that the $f_n$ form a Cauchy sequence
too. Since $L^1(M^\reg)$ is complete, there is a function $f\in
L^1(M^\reg)$ with $\lim_{n\to\infty}f_n=f$ and by continuity of
$\Psi$, we have $\Psi(f)=F$. Since $M^\reg$ (resp.
$G/H\times\Sigma^\reg$) is open and dense in $M$ (resp.
$G/H\times\Sigma$), we may identify $L^1(M^\reg)$ with $L^1(M)$
(resp. $L^1(G/H\times\Sigma^\reg)^W$ with
$L^1(G/H\times\Sigma)^W$). This proves the surjectivity of $\Psi$.
Using the open mapping theorem, we have thus established that
$\Psi$ is a continuous isomorphism between Banach spaces.

To (iii): Let us first assume that $f\in L^1(M)=L^1(M^\reg)$ has
support in some open subset $U\subset M^\reg$ such that
$\varphi^{-1}(U)$ is the union of $|W|$-many disjoint open subsets
$V_1,\dots, V_{|W|}$ of $G/H\times\Sigma^\reg$ each diffeomorphic
via $\varphi$ to $U$. The function $\Psi(f)$ has support in the
union of these $V_i$. By the classical transformation formula, we
then have
\begin{eqnarray*}
\int_M f(x)\, dx &=&\frac{1}{|W|}\int_{G/H\times\Sigma}f(g\cdot
s)|\det d\varphi(gH,s)|\, d(gH,s)\\
&\stackrel{(i)}{=}&\frac{1}{|W|}\int_\Sigma
\left(\int_{G/H}f(g\cdot s)\, d(gH)\right)\delta(s)\, ds.
\end{eqnarray*}
In the proof of (i), we have seen that $|\det
d\varphi(gH,s)|=|\det d\omega(eH)|=\delta(s)$, whence the last
equation above is true.

In the general case, with $f\in L^1(M^\reg)$ arbitrary, we cover
$M^\reg$ by a system $\ku$ of local trivial charts like $U$ above.
Take a partition of unity subordinate to $\ku$, say
$(U_i,p_i)_{i\in I}$ with $\sum_{i\in I} p_i=1$ and without loss
of generality $U_i\in\ku$. Then $f=\sum_{i\in I} f\cdot p_i$ and
\begin{eqnarray*}
\int_M f(x)\, dx &=&\sum_{i\in I}\int_M f(x)p_i(x)\, dx\\
&=&\frac{1}{|W|}\sum_{i\in I}\int_\Sigma \left(\int_{G/H}f(g\cdot
s)p_i(g\cdot s)\, d(gH)\right)\delta(s)\, ds\\
&=&\frac{1}{|W|}\int_\Sigma \left(\int_{G/H}f(g\cdot s)\,
d(gH)\right)\delta(s)\, ds.
\end{eqnarray*}

To (iv): Restricting $\Psi$ to the closed subspace
$L^1(M)^G\subset L^1(M)$, we obtain the designated mapping up to a
factor for the case $p=1$, because $\Psi(f)$ then does not depend
on the first variable anymore and the image of $\Psi$ is clearly
seen to be $L^1(\Sigma)^W$. The stated integral formula is an
immediate consequence of (iii), as well as the statement
concerning the isometries for the various $L^p$-spaces.
\end{proof}

\section{Examples}\label{s:examples}
The following two examples of polar actions and their corresponding
volume scaling functions $\delta$ are well known in the literature
(cf. \cite{Duistermaat&Kolk, Helgason2}). In this context, the
theorem above is called Weyl's Integration formula.

\begin{example}\label{e:example1}
Let $G$ be a compact connected Lie group $G$ and consider the action
of $G$ on itself via \emph{conjugation}; that is $\varphi:G\times
G\to G, (g,x)\mapsto gxg^{-1}$. Note that unless $G$ has trivial
center, this action is not effective. However, in view to the remark
(vi) above, this is not a problem.

A section is given by any maximal torus $T$ of $G$. It is pretty
straightforward to show that
$\delta(t)=|\det(\id-\Ad_t)|_{\gothg/\gotht}|$ for any $t\in T$. In
the presence of a \emph{root space decomposition}
$\gothg=\gotht\oplus\sum_{\alpha\in P} L_\alpha$, we can compute
$\delta$ more explicitly. For this purpose we recall some basic
facts we need:

Fix a maximal torus $T$ of $G$ and let $\gotht$ denote its Lie
algebra. A (infinitesimal) \emph{root} $\alpha$ is an element of
$(\gotht^{\C})^*$, the dual space of $\gotht^\C$, such that
$\gothg_\alpha:=\{Y\in\gothg^{\C}\mid \ad_X(Y)=\alpha(X)\cdot Y
\text{ for all } X\in\gotht\}\ne 0$. If $\alpha$ is a root then so
is $-\alpha$.  Let $P$ be a choice of \emph{positive roots}; that is
$0\notin P$ and for each $\alpha\in P$ we have $-\alpha\notin P$.
Now each $L_\alpha$ is defined by
$L_\alpha:=(\gothg_\alpha\oplus\gothg_{-\alpha})\cap\gothg$. Note
that since $\ad_X$ is a skew endomorphism, with respect to some
$\Ad$-invariant inner product, we see that each $\alpha$ takes
values in the imaginary numbers only. It is a fact that each nonzero
$\gothg_\alpha$ has complex dimension one, whereas the corresponding
$L_\alpha$ has real dimension two.

We then have that $\Ad_{\exp(X)}$ with $X\in\gotht$ acts on each
nonzero $\gothg_\alpha$ as the multiplication operator
$e^{\alpha(X)}\in U(1)$. Hence, the action of $\Ad_{\exp(X)}$ on
$L_\alpha$ is the rotation through the angle $\alpha(X)/i$. The root
space decomposition above is invariant under $(\id-\Ad_t)$. If we
write an arbitrary $t\in T$ as $t=\exp(X)$ for some $X\in\gotht$, we
obtain
$$\delta(t)=\prod_{\alpha\in P}\left[(1-\cos\frac{\alpha(X)}{i})^2+\sin\frac{\alpha(X)}{i}\right]=4^{|P|}\prod_{\alpha\in P}\sin^2\frac{\alpha(X)}{2i}.$$
\end{example}

\begin{example}\label{e:example2}
We may also consider the corresponding \emph{adjoint representation}
$\Ad$ of $G$ on $\gothg$. In this case, a section is given by any
maximal Abelian subspace $\gotht\subset\gothg$ of $\gothg$ which
then corresponds to a maximal torus of $G$. Here,
$\delta(X)=|\det(\ad_X)|$. Now $\ad_X$ acts on $\gothg_\alpha$ as
the multiplication operator $\alpha(X)$. Since $\alpha(X)$ assumes
purely imaginary values only, the action of $\ad_X$ on $L_\alpha$ is
a rotation through ninety degrees scaled by the factor
$\frac{\alpha(X)}{i}$. Hence,
$$\delta(X)=(-1)^{|P|}\prod_{\alpha\in
P}(\alpha(X))^2=\prod_{\alpha\in P}|\alpha(X)|^2.$$
\end{example}

The following two examples are generalizations of the previous ones
and they can, for instance, be found in \cite{Helgason2}(Ch. I,
\S5). However, the notion of a polar action is not used there
although the concept was already known at the time. Surprisingly,
the bibliography of the book even lists the article of Dadok, where
the polar representations on $\R^n$ are classified\cite{Dadok 1985}.

\begin{example}\label{e:example3}
Let $M=G/K$ be a Riemannian symmetric space with $G=\iso(M)^\circ$,
the identity component of the isometry group of $M$, and $K=G_p$,
the isotropy subgroup of $G$ of some point $p\in M$. For simplicity
reasons we will assume that $M$ is either of compact type or
noncompact type. It is then known that the action of $K$ on $M=G/K$
by left translation is \mbox{(hyper-)}polar and a section is given
by any maximal flat $A$ of $M$ through $p$. The situation here
indeed generalizes example \ref{e:example1}, because every connected
compact Lie group is a compact Riemannian symmetric space. In order
to compute $\delta(a)$, we start with some preliminaries. Let
$\gothg=\gothk\oplus\gothp$ be the Cartan decomposition of $\gothg$
with respect to $K$ and assume the usual $\Ad_K$ invariant inner
product on $\gothg$ derived from the Killing form; i.e. if $B$
denotes the Killing form of $\gothg$ then in the compact case, the
inner product is $-B$ and in the noncompact case it is
$-B(\cdot,\theta\cdot)$, where $\theta$ denotes the Cartan
involution. We identify $T_pM$ with $\gothp$ and, under this
identification, let $\gotha:=T_pA$. Then $\gotha\subset\gothp$ is a
maximal Abelian subspace of $\gothp$. If $\pi:G\to G/K$ denotes the
canonical projection, then $d\pi(e):\gothg\to\gothp$ is the
corresponding projection with respect to the Cartan decomposition
and the above identification. Let $L$ denote the centralizer of $A$
in $K$. Then $\omega_{aK}:K/L\to M,\ kL\mapsto kaK$ and we have
\begin{eqnarray*}
d\omega_{aK}(eL)(X+\gothl)&=&\left.\frac{d}{dt}\right|_{t=0}\omega_{aK}((\exp
tX)L)=\left.\frac{d}{dt}\right|_{t=0}(\exp tX)aK\\
&=&\left.\frac{d}{dt}\right|_{t=0}aa^{-1}(\exp tX)aK\\
&=&\left.\frac{d}{dt}\right|_{t=0}a(\exp t\Ad_{a^{-1}}X)K\\
&=&dl_a(p)\circ d\pi(e)(\Ad_{a^{-1}} X).
\end{eqnarray*}
Since $dl_a(p)$ is an isometry, we have
$\delta(aK)=|\det(d\pi(e)\circ\Ad_{a^{-1}})|$. For further
computations, we decompose
$\gothp=\gotha\oplus\sum_{\alpha\in\Sigma^+}\gothp_\alpha$ and
$\gothk=\gothl\oplus\sum_{\alpha\in\Sigma^+}\gothk_\alpha$. Here
$\Sigma$ denotes the set of restricted roots and $\Sigma^+$ a
choice of positive roots. Furthermore,
$\gothp_\alpha=\{Y\in\gothp\mid (\ad_X)^2 Y=\alpha(X)^2Y \text{
for all } X\in\gotha\}$ and $\gothk_\alpha=\{Y\in\gothk\mid
(\ad_X)^2 Y=\alpha(X)^2Y \text{ for all } X\in\gotha\}$ and we put
$m_\alpha:=\dim\gothp_\alpha=\dim\gothk_\alpha$; see e.g.
\cite{Helgason1} Ch. VII, \S 11 for the details. Now
$d\pi(e)\circ\Ad_{a^{-1}}$ leaves the decomposition of $\gothp$
above invariant and acts on each direct summand $\gothp_\alpha$ as
the operator $-\sinh\alpha(H)\cdot\frac{\ad_H}{\alpha(H)}$ in case
that $M$ is of noncompact type and it acts as the operator
$\sin\frac{\alpha(H)}{i}\cdot\frac{\ad_H}{\alpha(H)/i}$ in case
that $M$ is of compact type, where $a=\exp(H)$. In fact, if $M$ is
of noncompact type, then $\ad_H$ is a symmetric operator and thus
it has only real eigenvalues, whereas if $M$ is of compact type,
$\ad_H$ is skew-symmetric and has only purely imaginary
eigenvalues. In the noncompact case, we have for any
$X\in\gothp_\alpha$:
\begin{eqnarray*}
\Ad_{a^{-1}}X&=&e^{-\ad_H}(X)=\sum_{k=0}^\infty\frac{1}{k!}(-\ad_H)^k(X)\\
&=&-\sum_{k=0}^\infty\frac{1}{(2k+1)!}(\ad_H)^{2k+1}(X)+\sum_{k=0}^\infty\frac{1}{(2k)!}(\ad_H)^{2k}(X)\\
&=&-\frac{1}{\alpha(H)}\sum_{k=0}^\infty\frac{1}{(2k+1)!}\alpha(H)^{2k+1}\ad_H(X)+\sum_{k=0}^\infty\frac{1}{(2k)!}\alpha(H)^{2k}\cdot X\\
&=&-\sinh\alpha(H)\cdot\frac{\ad_H(X)}{\alpha(H)}+\cosh\alpha(H)\cdot
X.
\end{eqnarray*}
Clearly, after applying the projection $d\pi(e)$, only the first
term survives. Now $\frac{\ad_H}{\alpha(H)}$ is an isometry
between $\gothk_\alpha$ and $\gothp_\alpha$ and hence
$$\delta(aK)=\prod_{\alpha\in\Sigma^+}|\sinh\alpha(H)|^{m_\alpha}.$$
An analogue computation in the compact case reveals:
$$\delta(aK)=\prod_{\alpha\in\Sigma^+}\left|\sin\frac{\alpha(H)}{i}\right|^{m_\alpha}.$$
\end{example}

\begin{example}
As before, there is a Lie algebra version of the previous class of
examples, the so called \emph{s-representations}: The action of $K$
on $\gothp$ (see example \ref{e:example3} for the notation) given by
$k\cdot X:=\Ad_k(X)$ is \mbox{(hyper-)}polar. A section is any
Cartan subspace $\gotha$ of $\gothp$. Quite straightforward, we then
have $d\omega_H(eL)=\ad_H$ for any $H\in\gotha$. In view to the root
space decomposition in example \ref{e:example3}, we obtain
$$\delta(H)=\prod_{\alpha\in\Sigma^+}|\alpha(H)|^{m_\alpha}.$$
\end{example}

\begin{example}
The last family of examples we give here are the so called Hermann
actions (see \cite{Hermann}) which again generalize example
\ref{e:example3}. Assuming the notation of example \ref{e:example3}
with $G$ compact, let $H\subset G$ be another symmetric subgroup,
i.e. $G/H$ is a Riemannian symmetric space too. Then $H$ acts on
$M=G/K$ by left translation and this action is hyperpolar. A section
is given in the following way: In addition to the Cartan
decomposition $\gothg=\gothk\oplus\gothp$ with respect to $K$, we
also have the Cartan decomposition $\gothg=\gothh\oplus\gothm$ with
respect to $H$. Let $\gotha\subset\gothm\cap\gothp$ be a maximal
Abelian subspace. Then $A:=\exp(\gotha)$ is a flat section for the
action of $H$ on $M$ (details can be found in \cite{Conlon}). In
opposition to the other examples, the author does not know yet
wether it is possible to express the volume scaling function
$\delta$ in terms of roots in general, or not. A computation, as in
example \ref{e:example3}, however reveals that
$$\delta(aK)=|\det(d\pi(e)\circ\Ad_{a^{-1}})|,$$
where $a\in A$. The difficulty which arises is that
$d\pi(e)\circ\Ad_{a^{-1}}:\gothh/\gothz(A)\to\gothp$ and it is not
clear how to relate $\gothh/\gothz(A)$ with the normal space of
$\gotha$ in $\gothp$ in a more or less canonical way like we did
before, in the case when $\gothh=\gothk$.
\end{example}
%%%%%%%%%%%%%%%%%%%%%%%%%%%%%%%%%%%%%%%%%%%%%%%%%%%%%%%%%%%%%%%%%%%%%%%%%%%%%%%%

\end{document}